\documentclass{article}
\usepackage{amssymb}

\input{tcilatex}
\begin{document}

\textbf{\ \ \ \ \ \ \ \ \ \ \ \ Ap\'{e}ry, Bessel, Calabi-Yau and Verrill.}

\ \ \ \ \ \ \ \ \ \ \ \ \ \ \ \ \ \ \ \ Gert Almkvist

\textbf{Introduction.}

In [4] Bailey et al (among other things) study the Bessel moments

\[
c_{m,k}=\dint\limits_{0}^{\infty }x^{k}K_{0}(x)^{m}dx 
\]%
Here \ $K_{0}(x)$ \ is a certain Bessel function that conveniently can be
defined by%
\[
K_{0}(x)=\dint\limits_{0}^{\infty }e^{-x\cosh (t)}dt 
\]%
This leads to another representation (in Ising theory)%
\[
c_{m,k}=\frac{k!}{2^{m}}\dint\limits_{0}^{\infty
}...\dint\limits_{0}^{\infty }\frac{dx_{1}...dx_{m}}{(\cosh
(x_{1})+...+\cosh (x_{m}))^{k+1}} 
\]%
(historically it was the other way around).

In J.Borwein-Salvy [5] recursion formulas for the \ $c_{m,k}$ \ are derived (%
$m$ \ fixed). In the first section these recursions are studied in more
detail. E.g. if we define%
\[
d_{n}=\frac{16^{n}}{n!^{2}}c_{4,2n+1} 
\]%
we find an Ap\'{e}ry-like recursion (compare [3]) and recognize formulas
from [1] and [3]. Similar transformations of $c_{5,2n+1}$ \ lead to a 4-th
order differential equations whose mirror at \ $x=\infty $ \ is a Calabi-Yau
equation found by Verrill (\#34 in the "big table" [2]). This is also the
case with \ $c_{6,2n+1}$ \ where the differential equation at $\infty $ is
of order 5 (also found by Verrill) with a Calabi-Yau pullback of order 4
(\#130 in [2]).

There is an infinite sequence of differential equations of Verrill where the
coefficients are%
\[
A_{n}^{(m)}=\sum_{i_{1}+...i_{m}=n}(\frac{n!}{i_{1}!...i_{m}!})^{2} 
\]%
In [6] she gives a rather complicated formula for computing the recursion of
\ $A_{n}^{(m)}.$ In the second part we simplify this essentially using ideas
in J.Borwein-Salvy [5].

In the last section we prove the

\textbf{Main Theorem} \ For \ $m\geqslant 3$ \ we have

\[
y=\sum_{n=0}^{\infty }\frac{1}{4^{n}n!^{2}}c_{m,2n+1}x^{n} 
\]%
and%
\[
w=\sum_{n=0}^{\infty }A_{n}^{(m)}x^{-(n+1)} 
\]%
satisfy the same Picard-Fuchs differential equation of order \ $m_{+}=m/2$ \
if \ $m$ \ is even and \ $=(m+1)/2$ \ if \ $m$ \ is odd.This equation is
easily found by a Maple program.

There is a simplified version of this result for Bessel fans:

The differential equation satisfied by%
\[
y=\sum_{n=0}^{\infty }c_{m,2n}x^{2n} 
\]%
also has the solution%
\[
w=x^{-1}I_{0}(x^{-1})^{m} 
\]%
This depends on the identity%
\[
I_{0}(4\sqrt{x})^{m}=\sum_{i_{1}+...i_{m}=n}\frac{1}{i_{1}!^{2}...i_{m}!^{2}}%
x^{n} 
\]

\textbf{I. Some examples.}

\textbf{Four Bessel Functions}\bigskip

On p.13 in [4] Bailey et al define%
\[
c_{4,2n+1}=\dint\limits_{0}^{\infty }x^{2n+1}K_{0}(x)^{4}dx 
\]%
where \ $K_{0}$ is a Bessel function. In [5] the following recursion is
derived%
\[
64(k+3)c_{4,k+4}-4(k+2)(5k^{2}+20k+23)c_{4,k+2}+(k+1)^{5}c_{4,k}=0 
\]%
We make the substitution%
\[
d_{n}=\frac{16^{n}}{n!^{2}}c_{4,2n+1} 
\]%
and get the recursion%
\[
(n+2)^{3}d_{n+2}-2(2n+3)(5n^{2}+15n+12)d_{n+1}+64(n+1)^{3}d_{n}=0 
\]%
Then 
\[
y=\sum_{n=0}^{\infty }d_{n}x^{n} 
\]%
satisfies the differential equation where \ $\theta =x\frac{d}{dx}$ 
\[
\theta ^{3}-2x(2\theta +1)(5\theta ^{2}+5\theta +2)+64x^{2}(\theta +1)^{3} 
\]%
which we recognize as equation $(\alpha )$ \ in [1]. Then%
\[
A_{n}=\sum_{k=0}^{n}\binom{n}{k}^{2}\binom{2k}{k}\binom{2n-2k}{n-k} 
\]%
satisfies the recursion with initial values \ $A_{-1}=0$, $A_{0}=1.$ Let \ $%
B_{n}$ be the solution with \ $B_{0}=0,$ $B_{1}=1.$ Then we have

\textbf{Theorem. }We have 
\[
d_{n}=\frac{7}{8}A_{n}\zeta (3)-3B_{n} 
\]

\textbf{Proof. }In [4] we find \ $c_{4,1}=\frac{7}{8}\zeta (3)$ \ and \ $%
c_{4,3}=\frac{7}{32}\zeta (3)-\frac{3}{16}$ \ giving \ $d_{0}=\frac{7}{8}%
\zeta (3)$ \ and \ $d_{1}=\frac{7}{2}\zeta (3)-3$ \ Then we use the
recursion.

\bigskip

We want to find the asymptotic behaviour of \ $A_{n}$ and \ $d_{n}$ as $%
n\rightarrow \infty $. Making the Ansatz 
\[
A_{n}=Cn^{b}\lambda ^{n} 
\]%
in the recursion we find \ $\lambda =16$ \ or \ $\lambda =4$ \ and \ $b=-%
\frac{3}{2}\cdot $\ Numerical experiments suggest%
\[
A_{n}\sim 0.36\frac{16^{n}}{n^{3/2}} 
\]%
and%
\[
d_{n}\sim 0.7\frac{4^{n}}{n^{3/2}} 
\]%
This gives%
\[
\frac{7}{24}\zeta (3)-\frac{B_{n}}{A_{n}}\sim \frac{C}{4^{n}} 
\]%
which proves%
\[
\frac{B_{n}}{A_{n}}\rightarrow \frac{7}{24}\zeta (3) 
\]%
\[
\]

\bigskip \textbf{Remark.} The differential equation%
\[
\theta ^{3}-2x(2\theta +1)(5\theta ^{2}+5\theta +2)+64x^{2}(\theta +1)^{3} 
\]%
is self dual at infinity and the coefficients can be written \ (H.Verrill,
[6])%
\[
A_{n}=\sum_{i+j+k+l=n}(\frac{n!}{i!j!k!l!})^{2} 
\]%
\[
\]

\textbf{Five Bessel functions.}

Consider%
\[
c_{5,2n+1}=\dint\limits_{0}^{\infty }x^{2n+1}K_{0}(x)^{5}dx 
\]%
Then using the ideas of [5] we find the recursion%
\[
225c_{5,n+6}-(259n^{2}+1554n+2435)c_{5,n+4} 
\]%
\[
+(35n^{4}+280n^{3}+882n^{2}+1288n+731)c_{5,n+2}-(n+1)^{6}c_{5,n}=0 
\]%
Make the substitution%
\[
d_{n}=\frac{15^{2n}}{n!^{2}}c_{5,2n+1} 
\]%
which gives the recursion%
\[
n^{2}(n-1)^{2}d_{n}=4(n-1)^{2}(259n^{2}-518n+285)d_{n-1} 
\]%
\[
-3600(35n^{4}-210n^{3}+483n^{2}-504n+201)d_{n-2}+3240000(n-2)^{4}d_{n-3} 
\]%
\[
\]%
Let \ $A_{n}$ \ be the solution of the recursion with initial values \ $%
A_{0}=1,$ $A_{1}=0,$ $A_{2}=0.$ Similarly let \ $B_{n}$ and \ $C_{n}$ \ be
solutions with \ $B_{0}=0,$ $B_{1}=1,$ $B_{2}=0,$ \ $C_{0}=0,$ $C_{1}=0,$ $%
C_{2}=1$ \ respectively. Then%
\[
d_{n}=A_{n}s+225B_{n}t+C_{n}(6750-4500s+64125t) 
\]%
where \ $s=c_{5,1}$ \ and \ $t=c_{5,3}$ . We also use the conjectured value
of \ $c_{5,5}=\frac{8}{15}-\frac{16}{45}s+\frac{76}{15}t.$ Unfortunately we
still do not know the exact values of \ $s$ \ and \ $t.$ Maybe they are
related to the Ap\'{e}ry limits of \ $\frac{B_{n}}{A_{n}}$ \ and \ $\frac{%
C_{n}}{A_{n}}$

\textbf{A related Calabi-Yau equation.}

With \ $\theta =x\frac{d}{dx}$ the differential equation satisfied by%
\[
y=\sum_{n=0}^{\infty }d_{n}x^{n} 
\]%
is%
\[
\theta ^{2}(\theta -1)^{2}-4x\theta ^{2}(259\theta
^{2}+26)+3600x^{2}(35\theta ^{4}+70\theta ^{3}+63\theta ^{2}+28\theta
+5)-3240000x^{3}(\theta +1)^{4} 
\]%
The last factor cointains \ $(\theta +1)^{4}$ \ which suggests that
transforming the equation to \ $x=\infty $ \ could give a Calabi-Yau
equation. This is indeed the case: The substitutions \ $\theta
\longrightarrow -\theta -1$ \ and \ $x\longrightarrow 900x^{-1}$ \ give%
\[
\theta ^{4}-x(35\theta ^{4}+70\theta ^{3}+63\theta ^{2}+28\theta +5) 
\]%
\[
+x^{2}(\theta +1)^{2}(259\theta ^{2}+518\theta +285)-225x^{3}(\theta
+1)^{2}(\theta +2)^{2}, 
\]%
an equation found by Helena Verrill [6] . It has \#34 in the big table [2]
and has the analytic solution

\[
y=\sum_{n=0}^{\infty }a_{n}x^{n} 
\]%
where%
\[
a_{n}=\sum_{i+j+k+l+m=n}(\frac{n!}{i!j!k!l!m!})^{2} 
\]

\textbf{Six Bessel functions.}

Consider 
\[
c_{6,k}=\dint\limits_{0}^{\infty }x^{k}K_{0}(x)^{6}dx 
\]%
As above we have%
\[
2304(k+4)c_{6,k+6}-16(k+3)(49k^{2}+294k+500)c_{6,k+4} 
\]%
\[
+8(k+2)(7k^{4}+56k^{3}+182k^{2}+280k+171)c_{6,k+2}-(k+1)^{7}c_{6,k}=0 
\]%
With the substitution%
\[
d_{n}=\frac{48^{2n}}{n!^{2}}c_{6,2n+1} 
\]%
we have the recursion%
\[
(2n+5)(n+3)^{2}(n+2)^{2}d_{n+3}-32(n+2)^{3}(196n^{2}+784n+843)d_{n+2} 
\]%
\[
+64\cdot 48^{2}(2n+3)(14n^{4}+84n^{3}+196n^{2}+210n+87)d_{n+1}-128\cdot
48^{4}(n+1)^{5}d_{n}=0 
\]%
\[
\]%
Consider the three solutions \ $A_{n},B_{n},C_{n}$ \ with initial values%
\[
A_{0}=1,A_{1}=0,A_{2}=0 
\]%
\[
B_{0}=0,B_{1}=1,B_{2}=0 
\]%
\[
C_{0}=0,C_{1}=0,C_{2}=1 
\]%
\[
\]%
respectively. Let \ $c_{6,1}=s,$ $c_{6,3}=t$ . Then \ $c_{6,5}=\frac{5}{48}-%
\frac{1}{36}s+\frac{85}{72}t$\bigskip\ is conjectured. Then we have%
\[
d_{n}=A_{n}s+2304B_{n}t+C_{n}(138240-36864s+1566720t) 
\]

\textbf{A related Calabi-Yau equation.}

Let%
\[
y=\sum_{n=0}^{\infty }d_{n}x^{n} 
\]%
Then \ $y$ satisfies the differential equation%
\[
\theta ^{2}(\theta -1)^{2}(2\theta -1)-32x\theta ^{3}(196\theta ^{2}+59) 
\]%
\[
+64\cdot 48^{2}x^{2}(2\theta +1)(14\theta ^{4}+28\theta ^{3}+28\theta
^{2}+14\theta +3)-128\cdot 48^{4}x^{3}(\theta +1)^{5} 
\]%
\[
\]%
We find the mirror equation at \ $x=\infty $ \ via the substitution \ $%
\theta \longrightarrow -\theta -1$ \ and \ $x\longrightarrow 96^{2}x^{-1}$%
\[
\theta ^{5}-2x(2\theta +1)(14\theta ^{4}+28\theta ^{3}+28\theta
^{2}+14\theta +3) 
\]%
\[
+4x^{2}(\theta +1)^{3}(196\theta ^{2}+392\theta +255)-1152x^{3}(\theta
+1)^{2}(\theta +1)^{2}(2\theta +3) 
\]%
This we recognize as \#130 in the big table. It was found by H.Verrill [6].
The coefficients are

\[
A_{n}=\sum_{i+j+k+l+m+s=n}(\frac{n!}{i!j!k!l!m!s!})^{2} 
\]

\textbf{Seven Bessel functions.}

Let%
\[
d_{n}=\frac{105^{2n}}{n!^{2}}c_{7,2n+1} 
\]%
Then%
\[
y=\sum_{n=0}^{\infty }d_{n}x^{n} 
\]%
satisfies%
\[
\theta ^{2}(\theta -1)^{2}(\theta -2)^{2}-8x\theta ^{2}(\theta
-1)^{2}(6458\theta ^{2}-6458\theta +2589) 
\]%
\[
+48\cdot 105^{2}x^{2}\theta ^{2}(658\theta ^{4}+396\theta ^{2}+17) 
\]%
\[
-64\cdot 105^{4}x^{3}(84\theta ^{6}+252\theta ^{5}+378\theta ^{4}+336\theta
^{3}+180\theta ^{2}+54\theta +7) 
\]%
\[
+256\cdot 105^{6}x^{4}(\theta +1)^{6} 
\]

The transformation to \ infinity by \ $\theta \longrightarrow -\theta -1$ \
and \ $x\longrightarrow 210^{2}x^{-1}$ \ gives%
\[
\theta ^{6}-x(84\theta ^{6}+252\theta ^{5}+378\theta ^{4}+336\theta
^{3}+180\theta ^{2}+54\theta +7) 
\]%
\[
3x^{2}(\theta +1)^{2}(658\theta ^{4}+2632\theta ^{3}+4344\theta
^{2}+3424\theta +1071) 
\]%
\[
-2x^{3}(\theta +1)^{2}(\theta +2)^{2}(6458\theta ^{2}+19374\theta +15505) 
\]%
\[
+105^{2}x^{4}(\theta +1)^{2}(\theta +2)^{2}(\theta +3)^{2} 
\]%
with solution%
\[
y=\sum_{n=0}^{\infty }A_{n}x^{n} 
\]%
where%
\[
A_{n}=\sum_{i+j+k+l+m+p+s=n}(\frac{n!}{i!j!k!l!m!p!s!})^{2} 
\]%
\[
\]

\textbf{II. Sums of squares of generalized binomial coefficients.}

In [6] Verrill has given a rather complicated formula for the recursion of 
\[
A_{n}^{(k)}=\sum_{i_{1}+i_{2}+...+i_{k}=n}(\frac{n!}{i_{1}!i_{2}!...i_{k}!}%
)^{2}
\]%
We will instead consider 
\[
a_{n}^{(k)}=\frac{A_{n}}{n!^{2}}=\sum_{i_{1}+i_{2}+...i_{k}=n}\frac{1}{%
i_{1}!^{2}i_{2}!^{2}...i_{k}!^{2}}
\]%
Consider 
\[
y=\sum_{j=0}^{\infty }\frac{x^{j}}{j!^{2}}
\]%
Then $\ y$ \ satisfies the differential equation 
\[
\theta ^{2}-x
\]%
Actually 
\[
y(x)=I_{0}(4\sqrt{x})
\]%
Then%
\[
w=y^{m}=\sum_{n=0}^{\infty }a_{n}^{(m)}x^{n}
\]%
Using Lemma 3 in J.Borwein and Salvy [5] we find the following Maple program
for computing the differential equation for \ $w$ \ for all \ $m.$%
\[
\text{S:=proc(m) local M,k; M(0):=1; M(1):=t; for k to m do}
\]%
\[
\text{M[k+1]:=x*diff(M[k],x)+M[k]*t-k*(m-k+1)*x*M[k-2]; od;}
\]%
\[
\text{series(expand(M[m+1],x=0,infinity); end;}
\]%
\[
\]

Let \ $m_{+}=m/2$ \ if \ $m$ \ is even and \ $m_{+}=(m+1)/2$ \ if \ $m$ \ is
odd. Then write%
\[
S_{m}=\sum_{j=0}^{m_{+}}x^{j}Q_{j}(\theta ) 
\]%
Then the differential equation satisfied by%
\[
\sum_{n=0}^{\infty }A_{n}^{(m)}x^{n}=\sum_{n=0}^{\infty
}n!^{2}a_{n}^{(m)}x^{n} 
\]%
is given by%
\[
\theta ^{-2}\sum_{j=0}^{m_{+}}x^{j}\dprod\limits_{s=0}^{j-1}(\theta
+s)Q_{j}(\theta ) 
\]%
\[
\]

\textbf{III. Proof of the Main Theorem.}

The Bessel function \ $K_{0}(x)$ \ satisfies the differential equation \ $%
T_{m}(x,\theta )$ \ given by the Maple program%
\[
\]%
\[
\theta ^{2}-x^{2} 
\]%
Using Lemma 3 in Borwein-Salvy [5] again we obtain the differential equation
\ $T_{m}(x,\theta )$ \ satisfied by \ $K_{0}(x)^{m}$ \ given by the Maple
program

\[
\text{T:=proc(m) local L,k; L(0):=1; L(1):=t; for k to m do} 
\]%
\[
\text{L[k+1]:=x*diff(L[k],x)+L[k]*t-k*(m-k+1)*x*L[k-2]; od;} 
\]%
\[
\text{series(expand(L[m+1],x=0,infinity); end;} 
\]%
\[
\]

The crucial part of the proof is the following

\textbf{Lemma. }We have%
\[
M_{k}(x,\theta )=2^{-(k+1)}L_{k}(2\sqrt{x},2\theta ) 
\]

\textbf{Proof: }We use induction on \ $k$. Assume%
\[
M_{k-1}=2^{-k}L_{k-1}(2\sqrt{x},2\theta )\text{ \ and \ }%
M_{k}=2^{-(k+1)}L_{k}(2\sqrt{x},2\theta )
\]%
Then%
\[
M_{k+1}=x\frac{\partial M_{k}}{\partial x}+M_{k}\theta -xk(m-k+1)M_{k-1}
\]%
\[
=x2^{-(k+1)}\frac{\partial }{\partial x}L_{k}(2\sqrt{x},2\theta
)+2^{-(k+1)}L_{k}(2\sqrt{x},2\theta )\theta -x2^{-k}k(m-k+1)L_{k-1}(2\sqrt{x}%
,2\theta )
\]%
\[
=2^{-(k+1)}x\frac{1}{\sqrt{x}}\frac{\partial }{\partial (2\sqrt{x})}L_{k}(2%
\sqrt{x},2\theta )+2^{-(k+2)}L_{k}(2\sqrt{x},2\theta )2\theta -(2\sqrt{x}%
)^{2}2^{-(k+2)}k(m-k+1)L_{k-1}(2\sqrt{x},2\theta )
\]%
\[
=2^{-(k+2)}L_{k+1}(2\sqrt{x},2\theta )
\]%
\[
\]

The rest of the proof is merely book-keeping. Recall that \ 
\[
T_{m}(x,\theta )=\sum_{j=0}^{m_{+}}x^{2j}P_{j}(\theta ) 
\]%
annihilates \ $K_{0}(x)^{m}$. Then by the Maple program following Example 5
in [5] we find the recursion for \ $c_{m,k}$ \ by substituting \ $\theta
\longrightarrow -k-1-2j$ \ in \ $P_{j}(\theta ).$ Since \ $k=2n+1$ \ we get
\ $\theta \longrightarrow -2(n+1+j)$ \ Then with%
\[
d_{n}=\frac{1}{4^{n}n!^{2}}c_{m,2n+1} 
\]%
we get the following recursion for \ $d_{n}$%
\[
\sum_{n=0}^{m_{+}}n^{2}(n+1)^{2}...(n+j)^{2}4^{m_{+}-j}P_{j}(-2(n+1-j))N^{j}=0 
\]%
where \ $Nf(n)=f(n).$ \ Converting to the differential equation for \ $%
y=\sum d_{n}x^{n}$ \ we make the substitution \ $n\longrightarrow \theta -j$
\ and $N^{j}\longrightarrow x^{m_{+}-j}$ in the coefficient of \ $N^{j}$%
\[
\sum_{j=0}^{m_{+}}x^{m_{+}-j}\theta ^{2}(\theta -1)^{2}...(\theta
-j)^{2}4^{m_{+}-j}P_{j}(-2(\theta +1)) 
\]%
To get the differential equation at \ $\infty $ \ we make the substitution \ 
$\theta \longrightarrow -\theta -1$ \ and \ $x\longrightarrow x^{-1}$ \ and
we get%
\[
\sum_{j=0}^{m_{+}}x^{j}4^{j}\theta ^{2}(\theta +1)^{2}...(\theta
+j)^{2}P_{j}(2\theta )=\sum_{j=0}^{m_{+}}x^{j}\theta ^{2}(\theta
+1)^{2}...(\theta +j)^{2}Q_{j}(\theta ) 
\]%
which is the differential equation satisfied by%
\[
y=\sum_{n=0}^{\infty }A_{n}^{(m)}x^{n} 
\]

\textbf{Acknowledgements.}

I want to thank Wadim Zudilin who sent me the paper [4]. I also thank Jan
Gustavsson for doing some computations.

\textbf{References.}

\textbf{1.} G. Almkvist, W. Zudilin, Differential equations, mirror maps and
zeta values, in: Mirror Symmetry V, N. Yui, S.-T. Yau, and J. D. Lewis (eds), AMS/IP Stud. Adv. Math. 38 (International Press \& Amer. Math. Soc., Providence, RI, 2007) 481-515; arXiv: math/0402386

\textbf{2.} G. Almkvist, C. van Enckevort, D. van Straten, W. Zudilin, Tables of
Calabi-Yau equations, arXiv: math/0507430.

\textbf{3.} G. Almkvist, D. van Straten, W. Zudilin, Ap\'{e}ry limits of
differential equations of order 4 and 5, Banff 2006, to appear in Fields Comm. Publ., Vol~54.

\textbf{4.} D. H. Bailey. J. M. Borwein, D. Broadhurst, M. L. Glasser, Elliptic
integral evaluations of Bessel moments, arXiv: hep-th/0801089.

\textbf{5.} J. M. Borwein, B. Salvy, A proof of a recursion for Bessel moments,
inria-00152799

\textbf{6.} H. A. Verrill, Sums of squares of binomial coefficients, with
applications to Picard-Fuchs differential equations, math.CO/0407327.\ \ 

\bigskip

Institute of Algebraic Meditation

Fogdar\"{o}d 208

S-24333 H\"{o}\"{o}r

Sweden

gert.almkvist@yahoo.se

\end{document}